\documentclass{amsart}
\usepackage{bbm}
\usepackage{amsfonts}
\usepackage{amsmath, amsthm, amsfonts, amssymb}
\usepackage{mathrsfs}
\usepackage{amsmath, color}

\newtheorem{remark}{Remark}

\newtheorem{thm}{Theorem}[section]
\newtheorem{lem}{Lemma}[section]

\input{tcilatex}

\begin{document}
\title{boundedness of iterated spherical average on modulation spaces}
\author{Huang Qiang* and Fan Dashan }
\address[Huang Qiang]{Department of Mathematics, Zhejiang Normal University,
Jinhua 321000, China}
\address[Fan Dashan]{Department of Mathematics, Zhejiang Normal University,
Jinhua 321000, China}
\thanks{This work is supported by the NSF of China (Grant No.11801518) and
NSFZJ (~~No.LQ18A010005).\\
E-mail addresses: huangqiang0704@163.com(Q.Huang), fan@uwm.edu(D.Fan)}

\begin{abstract}
The spherical average $A_{1}(f)$ and its iteration $(A_{1})^{N}$ are
important operators in harmonic analysis and probability theory. Also $%
\Delta (A_{1})^{N}$ is used to study the $K$ functional in approximation
theory, where $\Delta $ is the Laplace operator. In this paper, we obtain
the sufficient and necessary conditions to ensure the boundedness of $\Delta
(A_{1})^{N}$ from the modulation space $M_{p_{1},q_{1}}^{s_{1}}$ to the
modulation space $M_{p_{2},q_{2}}^{s_{2}}$ for $1\leq
p_{1},p_{2},q_{1},q_{2}\leq \infty $ and $s_{1},s_{2}\in \mathbb{R}$.
\end{abstract}

\maketitle

\section{Introduction}

Let $\mathbb{S}^{n-1}$ be the unit sphere in the Euclidean space $%
%TCIMACRO{\U{211d} }%
%BeginExpansion
\mathbb{R}
%EndExpansion
^{n},$ $n\geq 2.$ We equip it with the normalized surface Lebesgue measure $%
d\sigma (y^{\prime })$. The average operator of functions $f$ on the unit
sphere is defined as
\begin{equation*}
A_{1}(f)(x)=\int_{\mathbb{S}^{n-1}}f(x-y^{\prime })d\sigma (y^{\prime }).
\end{equation*}

This operator has a profound background in harmonic analysis, dating back to
early 1970's (see \cite{Stein1},\cite{Stein2}). Moreover, it is closely
related to the study of random walks in high dimensional spaces,  which is
originated by Pearson \cite{Pearson} about 120 years ago. An $N$-steps
uniform walk in $\mathbb{R}^{n}$ starts at the origin and consists of $N$
independent steps of length 1, each of which is taken into a uniformly
random direction. It is known that the probability density function $p_{N}(%
\frac{n-2}{2},x)$ of such a random walk is the Fourier inverse of $%
(A_{1})^{N}$ (see \cite{Borwein}), where $(A_{1})^{N}$ denotes the $N$
iteration of $A_{1}.$

The operator $A_{1}$ also plays a significant role in the approximation
theory (see \cite{BDD}). Let $\Delta $ be the Laplacian. In order to obtain
some equivalent forms of the K-functional in $L^{p}(\mathbb{R}^{n})$ spaces,
Belinsky, Dai and Ditzian in \cite{BDD} study the iterates $(A_{1})^{N}$ for
positive integers $N$ and obtain the following theorem.

\textbf{Theorem A }(\cite{BDD}) Let $1\leq p\leq \infty $, $n\geq 2$ and $N>%
\frac{2(n+2)}{n-1}$.  The inequality
\begin{equation*}
\left\Vert \Delta (A_{1})^{N}(f)\right\Vert _{L^{p}(\mathbb{R}^{n})}\preceq
\Vert f\Vert _{L^{p}(\mathbb{R}^{n})}
\end{equation*}%
holds for all $f\in L^{p}(\mathbb{R}^{n})$.

Theorem A then raised the following question.

\textbf{Question 1 (\cite{BDD}):} Find the smallest positive integer $N$ to
guarantee the inequality
\begin{equation}
\left\Vert \Delta (A_{1})^{N}(f)\right\Vert _{L^{1}(\mathbb{R}^{n})}\preceq
\Vert f\Vert _{L^{1}(\mathbb{R}^{n})}.
\end{equation}

This question was addressed by Fan and Zhao in \cite{FZ} using the well
known estimates of wave operators (see \cite{Miyachi1}\cite{Peral}),  and
recently the question was completely solved by Fan, Lou and Wang in \cite%
{FLW} in  the following theorem.

\textbf{Theorem B }(\cite{FLW}). \textit{Let }$n\neq 3,5,$\textit{\ and }$N$%
\textit{\ be positive integers. The inequality  }%
\begin{equation*}
\left\Vert \Delta (A_{1})^{N}(f)\right\Vert _{L^{1}(\mathbb{R}^{n})}\preceq
\Vert f\Vert _{L^{1}(\mathbb{R}^{n})}
\end{equation*}%
\textit{holds if and only if }$N>\frac{n+3}{n-1}$\textit{.}

\textit{Let }$n=3,5,$ \textit{and }$N$\textit{\ be positive integers. The
inequality  }$\ $\textit{\ }%
\begin{equation*}
\left\Vert \Delta (A_{1})^{N}(f)\right\Vert _{L^{1}(\mathbb{R}^{n})}\preceq
\Vert f\Vert _{L^{1}(\mathbb{R}^{n})}
\end{equation*}%
\textit{holds if and only if }$N\geq \frac{n+3}{n-1}$\textit{.}

\bigskip

The aim of this article is to explore the behaves of $\Delta (A_{1})^{N}$ on
the modulation spaces $M_{p,q}^{s},$ where%
\begin{equation*}
(p,q,s)\in \lbrack 1,+\infty )\times \lbrack 1,+\infty )\times \mathbb{R}.
\end{equation*}%
We recall that the modulation space $M_{p,q}^{s}$ was introduced by
Feichtinger in \cite{F} and his initial aim was to measure smoothness of a
function or distribution in a way different from $L^{p}$ spaces. Nowadays,
spaces $M_{p,q}^{s}$ are recognized as a useful tool for studying functional
analysis and pseudo-differential operators (see \cite{B}\cite{CFS}\cite{Tj}%
). The original definition of the modulation space in \cite{F} is based on
the short-time Fourier transform and window function. In \cite{WH}, Wang and
Hudizk gave an equivalent definition of the discrete version on modulation
spaces by employing the frequency-uniform-decomposition. Later, people found
that the space $M_{p,q}^{s}$, with this discrete version, is a good working
frame to study boundedness of some operators and certain Cauchy problems of
nonlinear partial differential equations (see \cite{M}\cite{ZFG}\cite{HDC}%
\cite{RMW}). For example, the wave operator
\begin{equation*}
\widehat{e^{i|D|}f}=e^{i|\xi |}\widehat{f}
\end{equation*}%
is bounded in $L^{p}$ spaces if and only if $p=2$\ when $n\geq 2.$ However, $%
e^{i|D|}$ is bounded on modulation space $M_{p,q}^{s}$ for any  $p,q\in
\lbrack 1,+\infty )$ and $s\in \mathbb{R}$.

Motivated by these works, in this paper, we study boundedness of $\Delta
(A_{1})^{N}$ on modulation spaces and give the sufficient and necessary
conditions on the boundedness of $\Delta (A_{1})^{N}$ from $%
M_{p_{1},q_{1}}^{s_{1}}$ to $M_{p_{2},q_{2}}^{s_{2}}$ for $1\leq
p_{1},p_{2},q_{1},q_{2}\leq \infty $, $s_{1},s_{2}\in \mathbb{R}$. The
following theorem is our main result.

\begin{thm}
\label{thm1.1} \textit{Let }$\sigma =2-\frac{n-1}{2}N$\textit{\ and }$1\leq
p_{i},q_{i}\leq \infty $\textit{\ }$s_{i}\in R$\textit{\ for }$i=1,2$\textit{%
. When }$q_{1}\leq q_{2},$\textit{\ the iterated spherical average }$\Delta
(A_{1})^{N}$\textit{\ is bounded from }$M_{p_{1},q_{1}}^{s_{1}}(R^{n})$%
\textit{\ to }$M_{p_{2},q_{2}}^{s_{2}}(R^{n})$\textit{\ if and only if}%
\begin{equation*}
p_{1}\leq p_{2}\ \text{and}\quad s_{1}\geq s_{2}+\sigma .
\end{equation*}%
\textit{When }$q_{1}>q_{2},$\textit{\ the iterated spherical average }$%
\Delta (A_{1})^{N}$\textit{\ is bounded from }$M_{p_{1},q_{1}}^{s_{1}}(R^{n})
$\textit{\ to }$M_{p_{2},q_{2}}^{s_{2}}(R^{n})$\textit{\ if and only if}%
\begin{equation}
p_{1}\leq p_{2}\ \quad \text{and}\quad s_{1}+\frac{n}{q_{1}}>s_{2}+\sigma +%
\frac{n}{q_{2}}.
\end{equation}
\end{thm}

\begin{remark}
In above theorem, we can see that the smallest iterate step N which ensures $%
\Delta (A_{1})^{N}$ is bounded on modulation spaces $M_{p,q}^{s}(\mathbb{R}%
^{n})$ for all $(p,q,s)\in \lbrack 1,+\infty )\times \lbrack 1,+\infty
)\times \mathbb{R}$ is $\frac{4}{n-1}$, which is smaller than that in $L^{p}(%
\mathbb{R}^{n})$ spaces (see Theorem B). Moreover, our theorem finds the
sufficiency and necessity for the boundedness of $\Delta (A_{1})^{N}$ from $%
M_{p_{1},q_{1}}^{s_{1}}$ to $M_{p_{2},q_{2}}^{s_{2}}$ on full ranges of $%
1\leq p_{1},p_{2},q_{1},q_{2}\leq \infty $ and $s_{1},s_{2}\in \mathbb{R}$.
\end{remark}

\begin{remark}
When $n=1$, the average of sphere is reduced to
\begin{equation*}
a_{1}(f)(x)=\frac{1}{2}(f(x+1)+f(x-1)).
\end{equation*}%
Clearly, $a_{1}(f)$ and its iterates $a_{1}^{N}(f)$ are not in general
smoother than $f(x)$. However, with the increase of the dimension of space,
the spherical average $A_{1}(f)$ shares more regularity than $f(x)$.
Actually, our result also reflects this interesting phenomenon. If we choose
$n=1$ in Theorem \ref{thm1.1}, by the isomorphism property of modulation
spaces (see Proposition 2.1 ), $(A_{1})^{N}$ is bounded form $%
M_{p,q}^{s_{1}}(\mathbb{R}^{1})$ to $M_{p,q}^{s_{2}}(\mathbb{R}^{1})$ if and
only if $s_{2}\geq s_{1}$ for any iterate steps $N$. However, when $n\geq 2$
in Theorem \ref{thm1.1}, we can gain $\frac{n-1}{2}$ units of regularity in
each iterate step of $A_{1}$.
\end{remark}

The sufficiency part of the proof for Theorem 1.1 is somewhat routine with
the help of Bernstein's multiplier theorem. The necessity part of the proof
is quite involved. Based on the structure of $M_{p,q}^{s}$ and asymptotic
form of the Fourier transform of $\Delta (A_{1})^{N},$ we construct a
sequence of functions $\left\{ f_{k_{j},\lambda }\right\} $ to achieve the
necessary conditions.

This paper is organized as follows. In Section 2, we will introduce some
preliminary knowledge which includes some properties of modulation spaces
and some useful lemmas. The proof of Theorem \ref{thm1.1} will be presented
in Section 3.

Throughout this paper, we use the inequality $A\preceq B$ to mean that there
is a positive number $C$ independent of all main variables such that $A\leq
CB$, and use the notation $A\simeq B$ to mean $A\preceq B$ and $B\preceq A$.

\section{Preliminaries and Lemmas}

In this section, we give the definition and discuss some basic properties of
modulation spaces. Also, we will prove some estimates and lemmas which will
be used in our proof.

\textbf{Definition 2.1 }(Modulation spaces) Let $\varphi (\xi )$ be a smooth
function satisfying $\varphi (\xi )\equiv 1$ for $x\in \{x\in \mathbb{R}%
^{n}:|\xi |<\frac{1}{2}\}$, $\mathrm{supp}\varphi \subset \{\xi \in \mathbb{R%
}^{n}:|\xi |<\frac{3}{2}\}$ and $\{\varphi _{k}\}$ be a partition of the
unity satisfying the following conditions:
\begin{equation*}
\sum_{k\in Z^{n}}\varphi (\xi -k)=1,\varphi _{k}(\xi ):=\varphi (\xi -k)
\end{equation*}%
for any $\xi \in R^{n}$. And let
\begin{equation*}
\Box _{k}:=\mathcal{F}^{-1}\varphi _{k}\mathcal{F}.
\end{equation*}

With this frequency-uniform decomposition operator, we define the modulation
spaces \ $M_{p,q}^{s}(R^{n}),$ for $0<p,q\leq \infty ,$ $-\infty <s<\infty ,$
by
\begin{equation*}
M_{p,q}^{s}(R^{n}):=\left\{ f\in S^{\prime }:\Vert f\Vert
_{M_{p,q}^{s}(R^{n})}=\left( \sum_{k\in Z^{n}}\langle k\rangle ^{sq}\Vert
\Box _{k}f\Vert _{L^{p}(\mathbb{R}^{n})}^{q}\right) ^{\frac{1}{q}}<\infty
\right\} ,
\end{equation*}%
where $\langle k\rangle =\sqrt{1+|k|^{2}}$. See \cite{WH} for details.%
\newline
\newline
\textbf{Proposition 2.1 }(Isomorphism, see \cite{WH}) Let $0<p,q\leq \infty
,s,\tau \in \mathbb{R}$.
\begin{equation*}
J_{\sigma }=(I-\triangle )^{\frac{\tau }{2}}:M_{p,q}^{s}\rightarrow
M_{p,q}^{s-\tau }
\end{equation*}%
is an isomorphic mapping, where \ $I$ \ is the identity mapping and \ $%
\Delta $ \ is the Laplacian.\newline
\newline
\textbf{Proposition 2.2 }(Embedding, see \cite{WH}) For $0<p_{i},q_{i}\leq
\infty ,s_{i}\in \mathbb{R}$ ($i=1,2$), we have\newline
\begin{flalign}
M^{s_{1}}_{p_{1},q_{1}}\subset M^{s_{2}}_{p_{2},q_{2}}, \quad \text{if} \quad s_{1}\geq
s_{2},0<p_{1}\leq p_{2},0<q_{1}\leq q_{2}
\end{flalign}%
\begin{flalign}
\quad\quad\quad M^{s_{1}}_{p_{1},q_{1}}\subset M^{s_{2}}_{p_{2},q_{2}}, \quad \text{if} \quad s_{1}>s_{2},q_{1}>q_{2},s_{1}-s_{2}>n/q_{2}-n/q_{1}.
\end{flalign}

The Fourier multiplier is a linear operator $m(D)$ whose action on a test
function $f$ is formally defined by
\begin{equation*}
\widehat{m(D)f}(\xi )=m(\xi )\hat{f}(\xi ).
\end{equation*}%
The function $m(\xi )$ is called the symbol or multiplier of $m(D)$. Up to a
constant multiple, $m(D)$ is a convolution operator with the kernel
\begin{equation*}
K(x)=(m(\xi ))^{\vee }(x)=\int_{\mathbb{R}^{n}}m(\xi )e^{i\xi \cdot x}d\xi .
\end{equation*}%
By the Young inequality, we have
\begin{equation*}
\Vert m(D)f\Vert _{L^{p}}\preceq \left\Vert (m(\xi ))^{\vee }(x)\Vert
_{L^{1}}\right\Vert f\Vert _{L^{p}}
\end{equation*}%
for any $1\leq p\leq \infty $. We will use the following Bernstein
multiplier theorem to estimate $\Vert (m(\xi ))^{\vee }(x)\Vert _{L^{1}}$.

\begin{lem}
(Bernstein's multiplier theorem, see \cite{W1}) \label{lem2.1} Assume $%
0<p\leq 2$ and $\partial ^{\gamma }m(\xi )\in L^{2}$ for all multi-indices $%
\gamma $ with $|\gamma |\leq \left[ n(\frac{1}{p}-\frac{1}{2})\right] +1$.
We have
\begin{equation}
\left\Vert (m(\xi ))^{\vee }(x)\right\Vert _{L^{p}}\preceq
\sum\limits_{|\gamma |\leq \left[ n(\frac{1}{p}-\frac{1}{2})\right] +1}\Vert
\partial ^{\gamma }m(\xi )\Vert _{L^{2}}.
\end{equation}
\end{lem}

By checking the Fourier transform (see \cite{Stein2}), we have that
\begin{equation*}
\widehat{\Delta (A_{1})^{N}f}\simeq |\xi |^{2}(V_{\frac{n-2}{2}}(|\xi |))^{N}%
\widehat{f}(\xi )
\end{equation*}%
where
\begin{equation*}
V_{\delta }(r)=\frac{J_{\delta }(r)}{r^{\delta }}
\end{equation*}%
and $J_{\delta }(r)$ is the Bessel function of order $\delta $. Recall the
following asymptotic form of $J_{\delta }(r)$.

\begin{lem}
(\cite{Stein2})\label{lem2.2} Let $r>1$ and $\delta >-\frac{1}{2}$. For any
positive integer $L$ and $r\in \lbrack 1,\infty )$, we have
\begin{equation}
J_{\delta }(r)=\sqrt{\frac{2}{\pi r}}cos\left( r-\frac{\delta \pi }{2}-\frac{%
\pi }{4}\right) +\sum\limits_{j=1}\limits^{L}a_{j}e^{ir}r^{-\frac{1}{2}%
-j}+\sum\limits_{j=1}\limits^{L}b_{j}e^{-ir}r^{-\frac{1}{2}-j}+E(r)
\end{equation}%
where $a_{j}$ and $b_{j}$ are constants for all $j$, and $E(r)$ is a $%
C^{\infty }$ function satisfying
\begin{equation*}
\left\vert E^{(k)}(r)\right\vert \preceq r^{-\frac{1}{2}-L-1}
\end{equation*}%
for any $k=0,1,2...$
\end{lem}

\section{Proof of Theorem 1.1}

We start with showing the sufficiency of Theorem 1.1. By the definition of
modulation spaces, we need to estimate $\left\Vert \square _{k}\Delta
(A_{1})^{N}f\right\Vert _{L^{p_{2}}(\mathbb{R}^{n})}$ and obtain the
following lemma.

\begin{lem}
\label{lem3.1} Let $1\leq p_{2}\leq \infty $ and $\sigma =2-\frac{n-1}{2}N$.
Then
\begin{equation*}
\left\Vert \square _{k}\Delta (A_{1})^{N}f\right\Vert _{L^{p_{2}}(\mathbb{R}%
^{n})}\preceq \langle k\rangle ^{\sigma }\Vert \square _{k}f\Vert
_{L^{p_{2}}(\mathbb{R}^{n})}.
\end{equation*}
\end{lem}

\textbf{Proof:} For $\forall k\in \mathbb{Z}^{n}$, $\square _{k}\Delta
(A_{1})^{N}f$ is a convolution operator $m_{k}(D)(f)=\Omega _{k}(x)\ast f$,
where
\begin{equation}\label{eqo3.7}
\Omega _{k}(x)=\int_{\mathbb{R}^{n}}\varphi _{k}(\xi )|\xi |^{2}\left( V_{%
\frac{n-2}{2}}(|\xi |)\right) ^{N}e^{i\xi \cdot x}d\xi .
\end{equation}%
By the almost orthogonality of unit decomposition, there exists an integer $%
k_{0}(n)$ which depends only on $n$ such that $\varphi _{l}(\xi )\varphi
_{k}(\xi )=0$ when $|l-k|\geq k_{0}(n)$. Since
\begin{equation*}
\sum_{k\in Z^{n}}\square _{k}=I
\end{equation*}%
where $I$ is the identity operator, Young's inequality and Minkowski's
inequality yield
\begin{eqnarray*}
\left\Vert \square _{k}\Delta (A_{1})^{N}f\right\Vert _{L^{p_{2}}} &\leq
&\sum\limits_{l\in \mathbb{Z}^{n},|l-k|\leq k_{0}(n)}\left\Vert \square
_{l}\Delta (A_{1})^{N}\square _{k}f\right\Vert _{L^{p_{2}}} \\
&\preceq &\sum\limits_{l\in \mathbb{Z}^{n},|l-k|\leq k_{0}(n)}\left\Vert
\left( \varphi _{l}(\xi )|\xi |^{2}(V_{\frac{n-2}{2}}(|\xi |))^{N}\right)
^{\vee }\right\Vert _{L^{1}}\Vert \square _{k}f\Vert _{L^{p_{2}}}.
\end{eqnarray*}%
So, it suffices to estimate
\begin{equation*}
\sum\limits_{l\in \mathbb{Z}^{n},|l-k|\leq k_{0}(n)}\left\Vert \left(
\varphi _{l}(\xi )|\xi |^{2}(V_{\frac{n-2}{2}}(|\xi |))^{N}\right) ^{\vee
}\right\Vert _{L^{1}}
\end{equation*}%
for every $k\in \mathbb{Z}^{n}.$ Notice that the cardinality of
\begin{equation*}
\Lambda _{k}:=\{l\in \mathbb{Z}^{n}:|l-k|\leq k_{0}(n)\}
\end{equation*}%
is uniformly finite for all $k\in \mathbb{Z}^{n}$, and $\langle l\rangle
\simeq \langle k\rangle $ when $l\in \Lambda _{k}$. Therefore, we only need
to estimate the $L^{1}$ norm
\begin{equation*}
\left\Vert \left( \varphi _{l}(\xi )|\xi |^{2}(V_{\frac{n-2}{2}}(|\xi
|))^{N}\right) ^{\vee }\right\Vert _{L^{1}}=\|\Omega_{l}(x)\|_{L^{1}}
\end{equation*}
for $\langle l\rangle \simeq \langle k\rangle $.

When $|k|<100$, by the well known formula (\cite{Stein2})
\begin{equation*}
V_{\delta }(t)=O(1)\quad \text{if}\quad |t|<100
\end{equation*}%
we have that $|\Omega _{l}(x)|\preceq 1$ for $|x|<100$.

On the other hand, when $|x|\geq 100$, without loss of generality, we may
assume $|x_{1}|\geq \frac{|x|}{n}$. By the derivative formula of $V_{\delta
}(t)$
\begin{equation}
\frac{dV_{\delta }(t)}{dt}=-tV_{\delta +1}(t)
\end{equation}%
and taking integration by part on $\xi _{1}$ variable in (\ref{eqo3.7}), we obtain that
\begin{equation*}
|\Omega _{l}(x)|\preceq \frac{1}{|x_{1}|^{n+1}}\preceq \frac{1}{|x|^{n+1}}
\end{equation*}%
for $|x|\geq 100$. This estimate implies that $\Vert \Omega _{l}(x)\Vert
_{L^{1}}\preceq 1$ when $|k|<100$, since $\langle l\rangle \simeq \langle
k\rangle $.

Next, we study the case $|k|\geq100$. Choosing $L=1$ in Lemma \ref{lem2.2},
we have the following asymptotic form of $V_{\delta}(r)$
\begin{equation}
V_{\delta}(r)=r^{-\delta-\frac{1}{2}}\left(\sqrt{\frac{2}{\pi}}\cos\left(r-%
\frac{\delta\pi}{2}-\frac{\pi}{4}\right)  \right)+O(r^{-\delta-\frac{3}{2}})
\end{equation}
for $|r|>1$.

Therefore, when $|k|>100$ and $\langle l\rangle \simeq \langle k\rangle $,
we have
\begin{equation}
\left\vert V_{\delta }(|\xi |)^{N}\right\vert \preceq |l|^{(-\delta -\frac{1%
}{2})N}\simeq \langle k\rangle ^{(-\delta -\frac{1}{2})N}
\end{equation}%
for $\xi \in \mathrm{supp}\varphi _{l}(\xi )$. Now, by the chain rule and
the derivative formula of $V_{\delta }(t)$ , we obtain
\begin{equation}
\begin{split}
\frac{\partial }{\partial \xi _{i}}\left( V_{\delta }(|\xi |)\right) ^{N}&
=-\left( V_{\delta }(|\xi |)\right) ^{N-1}|\xi |\cdot V_{\delta +1}(|\xi
|)\cdot \frac{\xi _{i}}{|\xi |} \\
& =-\left( V_{\delta }(|\xi |)\right) ^{N-1}V_{\delta +1}(|\xi |)\cdot \xi
_{i}.
\end{split}
\label{eqn2}
\end{equation}%
By the asymptotic form of $V_{\delta }(r)$, we obtain that
\begin{equation}
\left\vert \frac{\partial }{\partial \xi _{i}}\left( V_{\delta }(|\xi
|)\right) ^{N}\right\vert \preceq |\xi |^{(-\delta -\frac{1}{2})(N-1)}|\xi
|^{-\delta -\frac{3}{2}}|\xi |\preceq |l|^{(-\delta -\frac{1}{2})N}\simeq
\langle k\rangle ^{(-\delta -\frac{1}{2})N}
\end{equation}%
for $\xi \in \mathrm{supp}\varphi _{l}(\xi )$.

Thus, $V_{\delta }(|\xi |)^{N}$ and $\frac{\partial }{\partial \xi _{i}}%
\left( V_{\delta }(|\xi |)\right) ^{N}$ share the same upper bound which is
$\langle l\rangle ^{(-\delta -\frac{1}{2})N}$, for any $%
\delta $ and $\xi \in \mathrm{supp}\varphi _{l}(\xi )$. By Lemma \ref{lem2.1}
(Bernstein's multiplier theorem) and the fact $\partial ^{\gamma }(|\xi
|^{2})\preceq |\xi |^{2-|\gamma |}$ for $|\gamma |\leq 2$ and $\partial
^{\gamma }(|\xi |^{2})=0$ for $|\gamma |>2$, we have that
\begin{eqnarray*}
\Vert \Omega _{l}(x)\Vert _{L^{1}} &=&\left\Vert \left( \varphi _{l}(\xi
)|\xi |^{2}\left( V_{\frac{n-2}{2}}(|\xi |)\right) ^{N}\right) ^{\vee
}\right\Vert _{L^{1}} \\
&\preceq &\sum\limits_{|\gamma |\leq \left[ \frac{n}{2}\right] +1}\left\Vert
\partial ^{\gamma }\left( \varphi _{l}(\xi )|\xi |^{2}(V_{\frac{n-2}{2}%
}(|\xi |))^{N}\right) \right\Vert _{L^{2}} \\
&\preceq &\sum\limits_{|\gamma |\leq \left[ \frac{n}{2}\right]
+1}\sum\limits_{\gamma _{1}+\gamma _{2}+\gamma _{3}=\gamma }\left\Vert
\partial ^{\gamma _{1}}\varphi _{l}(\xi )\cdot \partial ^{\gamma _{2}}|\xi
|^{2}\cdot \partial ^{\gamma _{3}}\left( V_{\frac{n-2}{2}}(|\xi |)\right)
^{N}\right\Vert _{L^{2}(\mathrm{supp}\varphi _{l}(\xi ))} \\
&\preceq &\sum\limits_{|\gamma |\leq \left[ \frac{n}{2}\right]
+1}\sum\limits_{\gamma _{1}+\gamma _{2}=\gamma }\left\Vert |\xi |^{2}\cdot
\partial ^{\gamma _{1}}\varphi _{l}(\xi )\cdot \partial ^{\gamma _{2}}\left(
V_{\frac{n-2}{2}}(|\xi |)\right) ^{N}\right\Vert _{L^{2}(\mathrm{supp}%
\varphi _{l}(\xi ))} \\
&\preceq &|l|^{2-\frac{n-1}{2}N}\simeq \langle k\rangle ^{2-\frac{n-1}{2}N}
\end{eqnarray*}%
for $l\in\Lambda_{k}$.

Combining all above estimates, by the definition of modulation spaces, we
have that
\begin{eqnarray*}
\left\Vert \triangle (A_{1})^{N}f\right\Vert _{M_{p_{2},q_{2}}^{s_{2}}}
&=&\left( \sum\limits_{k\in \mathbb{Z}^{n}}\langle k\rangle
^{s_{2}q_{2}}\left\Vert \Box _{k}\triangle (A_{1})^{N}f\right\Vert
_{L^{p_{2}}}^{q_{2}}\right) ^{\frac{1}{q_{2}}} \\
&=&\left( \sum\limits_{|k|<100}\langle k\rangle ^{s_{2}q_{2}}\left\Vert \Box
_{k}\triangle (A_{1})^{N}f\right\Vert
_{L^{p_{2}}}^{q_{2}}+\sum\limits_{|k|\geq 100}\langle k\rangle
^{s_{2}q_{2}}\left\Vert \Box _{k}\triangle (A_{1})^{N}f\right\Vert
_{L^{p_{2}}}^{q_{2}}\right) ^{\frac{1}{q_{2}}} \\
&\preceq &\left( \sum\limits_{|k|<100}\langle k\rangle
^{s_{2}q_{2}}\left\Vert \Box _{k}\triangle (A_{1})^{N}f\right\Vert
_{L^{p_{2}}}^{q_{2}}\right) ^{\frac{1}{q_{2}}}+\left( \sum\limits_{|k|\geq
100}\langle k\rangle ^{s_{2}q_{2}}\left\Vert \Box _{k}\triangle
(A_{1})^{N}f\right\Vert _{L^{p_{2}}}^{q_{2}}\right) ^{\frac{1}{q_{2}}} \\
&\preceq &\left( \sum\limits_{|k|<100}\langle k\rangle ^{(s_{2}+2-\frac{n-1}{%
2}N)q_{2}}\Vert \Box _{k}f\Vert _{L^{p_{2}}}^{q_{2}}\right) ^{\frac{1}{q_{2}}%
}+\left( \sum\limits_{|k|\geq 100}\langle k\rangle ^{(s_{2}+2-\frac{n-1}{2}%
N)q_{2}}\Vert \Box _{k}f\Vert _{L^{p_{2}}}^{q_{2}}\right) ^{\frac{1}{q_{2}}}
\\
&\preceq &\Vert f\Vert _{M_{p_{2},q_{2}}^{s_{2}+2-\frac{n-1}{2}N}}\ .
\end{eqnarray*}%
By the embedding properties of modulation spaces (Proposition 2.2), we can
easily obtain that
\begin{equation}
\Vert \triangle (A_{1})^{N}f\Vert _{M_{p_{2},q_{2}}^{s_{2}}}\preceq \Vert
f\Vert _{M_{p_{2},q_{2}}^{s_{2}+\sigma }}\preceq \Vert f\Vert
_{M_{p_{1},q_{1}}^{s_{1}}}
\end{equation}%
when%
\begin{equation*}
p_{1}\leq p_{2}\quad \text{and }\quad s_{1}\geq s_{2}+\sigma \quad\text{ if }\quad q_{1}\leq
q_{2}
\end{equation*}%
or
\begin{equation*}
p_{1}\leq p_{2}\quad \text{and}\quad s_{1}+\frac{n}{q_{1}}>s_{2}+\sigma +
\frac{n}{q_{2}}\quad\text{ if }\quad q_{1}>q_{2}.
\end{equation*}%
The sufficiency of Theorem \ref{thm1.1} is proved.

Turn to prove the necessity part of Theorem \ref{thm1.1}. We need the
following lemma.

\begin{lem}
\label{lem3.2} Let $1\leq p\leq \infty $. These exists a constant $\rho
=\rho (n)>0$ which depends only on $n$ and a subsequence $\{k_{j}\}\subseteq
\mathbb{Z}^{n}$ such that
\begin{equation}
\left\Vert \Box _{k_{j}}\Delta (A_{1})^{N}g_{k_{j}}\right\Vert
_{L^{p}}\simeq \langle k_{j}\rangle ^{2-\frac{n-1}{2}N}\Vert g_{k_{j}}\Vert
_{L^{p}},
\end{equation}%
where $\{g_{k_{j}}(x)\}$ is a sequence of Schwartz function with $\mathrm{%
supp}\widehat{g_{k_{j}}}(\xi )\subset \{\xi \in \mathbb{R}^{n}:|\xi
-k_{j}|\leq \rho \}$.
\end{lem}

\textbf{Proof:} By the same method as in Lemma \ref{lem3.1}, it is easy to
get
\begin{equation}
\left\Vert \Box _{k}\Delta (A_{1})^{N}g_{k}\right\Vert _{L^{p}}\preceq \Vert
\Omega _{k}(x)\Vert _{L^{1}}\Vert g_{k}\Vert _{L^{p}}\preceq \langle
k\rangle ^{2-\frac{n-1}{2}N}\Vert g_{k}\Vert _{L^{p}}  \label{eqo3.1}
\end{equation}%
for all $k\in \mathbb{Z}^{n}$. Thus, we only need to prove the inverse
inequality. By Lemma \ref{lem2.2}, we have
\begin{eqnarray*}
V_{\frac{n-2}{2}}(r) &=&r^{-\frac{n-1}{2}}\sqrt{\frac{2}{\pi }}\cos (r-\frac{%
n\pi }{4}+\frac{\pi }{4})+O(r^{-\frac{n+1}{2}}) \\
&=&r^{-\frac{n-1}{2}}\sqrt{\frac{2}{\pi }}\sin (r-\frac{n\pi }{4}+\frac{3\pi
}{4})+O(r^{-\frac{n+1}{2}})
\end{eqnarray*}%
for $|r|>1$. We consider
$$u(r):=\sin (r-\frac{n\pi }{4}+\frac{3\pi }{4})$$
in every semiperiod $r-\frac{n\pi }{4}+\frac{3\pi }{4}\in \lbrack j\pi
,(j+1)\pi ]$, $(j=1,2,...).$

Choosing $\varepsilon _{0}=\sin (0.07)$, we have
\begin{equation}
\left\vert u(r)\right\vert \geq \varepsilon _{0}
\end{equation}%
for $r-\frac{n\pi }{4}+\frac{3\pi }{4}\in \lbrack j\pi +0.07,j\pi +\pi -0.07]
$, which is equivalent to
\begin{equation*}
r\in \left[ j\pi +\frac{n\pi }{4}-\frac{3\pi }{4}+0.07,\text{ \ }j\pi +\pi +%
\frac{n\pi }{4}-\frac{3\pi }{4}-0.07\right] .
\end{equation*}%
By Lemma \ref{lem3.3} (the lemma will be proved later), for every $j\in
\mathbb{N}^{+}$, the set
\begin{equation*}
\Lambda _{1,j}:=\left\{ k\in \mathbb{Z}^{n}:|k|\in \left[ j\pi +\frac{n\pi }{%
4}-\frac{3\pi }{4}+0.07,\text{ }j\pi +\pi +\frac{n\pi }{4}-\frac{3\pi }{4}%
-0.07\right] \right\}
\end{equation*}%
is not empty. So, there exists a subsequence of integer $\{k_{j}\}$, such
that $k_{j}\in \Lambda _{1,j}$ and
\begin{equation*}
\left\vert u(|k_{j}|)\right\vert \geq \varepsilon _{0}.
\end{equation*}%
Moreover,
\begin{equation*}
\left\vert u^{\prime }(r)\right\vert =\left\vert \cos (r-\frac{n\pi }{4}+%
\frac{3\pi }{4})\right\vert \leq 1,
\end{equation*}%
which means that
\begin{equation*}
|u(r)|\geq \frac{\varepsilon _{0}}{2}
\end{equation*}%
for
\begin{equation*}
r\in \left[ |k_{j}|-\frac{\varepsilon _{0}}{4},|k_{j}|+\frac{\varepsilon _{0}%
}{4}\right]
\end{equation*}%
and
\begin{equation*}
|k_{j}|\in \left[ j\pi +\frac{n\pi }{4}-\frac{3\pi }{4}+0.07,j\pi +\pi +%
\frac{n\pi }{4}-\frac{3\pi }{4}-0.07\right] .
\end{equation*}

%  By lemma XXX, for every $j\in\mathbb{N}^{+}$, the set
% $$\Lambda_{1,j}:=\left\{k\in\mathbb{Z}^{n}:|k|\in\left[j\pi+\frac{n\pi}{4}-\frac{3\pi}{4}+0.07,
% j\pi+\pi+\frac{n\pi}{4}-\frac{3\pi}{4}-0.07\right]\right\}$$
% is not empty. So, there exit a subsequence of integer $\{k_{j}\}$, such that $k_{j}\in\Lambda_{1,j}$ for all $j\in\mathbb{N}^{+}$
% and $|u(|r|)|\geq\frac{\varepsilon_{0}}{2}.$ for $r\in[|k_{j}|-\frac{\varepsilon_{0}}{4},|k_{j}|+\frac{\varepsilon_{0}}{4}]$

%Notice that the length of $$\left[2j\pi+\frac{n\pi}{4}-\frac{1}{2},2j\pi+\frac{n\pi}{4}+\frac{1}{2}\right]$$ is 1. So, there
%exit a positive integer $$n_{j}\in\left[2j\pi+\frac{n\pi}{4}-\frac{1}{2},2j\pi+\frac{n\pi}{4}+\frac{1}{2}\right].$$
%such that $|u(n_{j})|\geq\varepsilon_{0}$.

For the remainder $O(r^{-\frac{n+1}{2}})$ in the expansion of $V_{\frac{n-2}{%
2}}(r)$, it is obvious that, when $r$ is large enough,
\begin{equation*}
O(r^{-\frac{n+1}{2}})\leq \frac{\varepsilon _{0}}{4}r^{-\frac{n-1}{2}}.
\end{equation*}

Let $\epsilon ,\rho =\frac{\varepsilon _{0}}{4}$. We obtain that there
exist some constants $\epsilon ,\rho >0$ and a subsequence $%
\{k_{j}\}\subseteq \mathbb{Z}^{+}$ such that
\begin{equation}
|V_{\frac{n-2}{2}}(|\xi |)|\geq \epsilon |\xi |^{-\frac{n-1}{2}}
\label{eqo3.2}
\end{equation}%
for $\xi \in \{\xi :|\xi -k_{j}|\leq \rho \}$ when $j$ is large enough.
Moreover the subsequence $\{k_{j}\}\subseteq \mathbb{Z}^{n}$ satisfies
\begin{equation}
|k_{j}|\in \left[ j\pi +\frac{n\pi }{4}-\frac{3\pi }{4}+0.07,\text{ }j\pi
+\pi +\frac{n\pi }{4}-\frac{3\pi }{4}-0.07\right] ,  \label{eqo3.3}
\end{equation}%
when the positive integer $j$ is large enough.

Therefore, when $\xi \in \{\xi :|\xi -k_{j}|\leq \rho \}$ and $N\in \mathbb{Z%
}^{+}$, we have
\begin{equation*}
|V_{\frac{n-2}{2}}(|\xi |)^{-N}|\preceq |\xi |^{(\frac{n-1}{2})N}\simeq
\langle k_{j}\rangle ^{(\frac{n-1}{2})N}.
\end{equation*}%
Using the chain rule and the derivative formula of $V_{\delta }(t)$,
\begin{equation*}
\frac{\partial }{\partial \xi _{i}}\left( V_{\frac{n-2}{2}}(|\xi |)\right)
^{-N}=-\left( V_{\frac{n-2}{2}}(|\xi |)\right) ^{-(N+1)}V_{\frac{n-2}{2}%
+1}(|\xi |)\cdot \xi _{i}.
\end{equation*}%
By the asymptotic form of $V_{\delta }(r)$ and (\ref{eqo3.2}), we have
\begin{equation*}
\left\vert \frac{\partial }{\partial \xi _{i}}\left( V_{\frac{n-2}{2}}(|\xi
|)\right) ^{-N}\right\vert \preceq |\xi |^{(-\frac{n-1}{2})(-N-1)}|\xi |^{-%
\frac{n+1}{2}}|\xi |\simeq \langle k_{j}\rangle ^{(\frac{n-1}{2})N}.
\end{equation*}%
As a result, $V_{\frac{n-2}{2}}(|\xi |)^{-N}$ and $\frac{\partial }{\partial
\xi _{i}}\left( V_{\frac{n-2}{2}}(|\xi |)\right) ^{-N}$ share the same upper
bound which is $\langle k_{j}\rangle ^{(\frac{n-1}{2})N}$, for $\xi \in \{\xi :|\xi -k_{j}|\leq \rho \}.$

Let $\eta (\xi )$ be a smooth function with $\mathrm{supp}~{\eta _{k}(\xi )}%
\subset \{\xi :|\xi |\leq 2\rho \}$ and $\eta _{k}(\xi )\equiv 1$ for $\xi
\in \{\xi :|\xi |\leq \rho \}$. We define
\begin{equation*}
\eta _{k}(\xi )=\eta (\xi -k),k\in \mathbb{Z}^{n}.
\end{equation*}%
Notice that $\rho =\frac{\sin (0.07)}{4}<\frac{1}{4}$. Moreover, for the partition of the unity $\{\varphi_{k}\}$(see Definition 2.1), we have that
\begin{equation*}
\varphi _{k_{j}}(\xi )\widehat{g_{k_{j}}}(\xi )=\widehat{g_{k_{j}}}(\xi ),
\end{equation*}%
with
\begin{equation*}
\mathrm{supp}~\widehat{g_{k_{j}}}(\xi )\subset \{\xi :|\xi -k_{j}|\leq \rho
\}.
\end{equation*}%
Therefore, by the Bernstein multiplier theorem ( Lemma \ref{lem2.1}) and (%
\ref{eqo3.2}), we obtain that
\begin{eqnarray*}
\Vert g_{k_{j}}\Vert _{L^{p}} &=&\Vert (\widehat{g_{k_{j}}})^{\vee }\Vert
_{L^{p}} \\
&=&\left\Vert \left( \eta _{k_{j}}(\xi )\varphi _{k_{j}}(\xi )|\xi |^{-2}(V_{%
\frac{n-2}{2}}(|\xi |))^{-N}\cdot |\xi |^{2}(V_{\frac{n-2}{2}}(|\xi |))^{N}%
\widehat{g_{k_{j}}}(\xi )\right) ^{\vee }\right\Vert _{L^{p}} \\
&\preceq &\left\Vert (\eta _{k_{j}}(\xi )|\xi |^{-2}(V_{\frac{n-2}{2}}(|\xi
|))^{-N})^{\vee }\right\Vert _{L^{1}}\left\Vert \Box _{k_{j}}\Delta
(A_{1})^{N}g_{k_{j}}\right\Vert _{L^{p}} \\
&\preceq &\sum\limits_{|\gamma |\leq \left[ n(\frac{1}{p}-\frac{1}{2})\right]
+1}\left\Vert \partial ^{\gamma }(\eta _{k_{j}}(\xi )|\xi |^{-2}(V_{\frac{n-2%
}{2}}(|\xi |))^{-N})\right\Vert _{L^{2}}\left\Vert \Box _{k_{j}}\Delta
(A_{1})^{N}g_{k_{j}}\right\Vert _{L^{p}} \\
&\preceq &\sum\limits_{|\gamma |\leq \left[ n(\frac{1}{p}-\frac{1}{2})\right]
+1}\sum\limits_{\gamma _{1}+\gamma _{2}+\gamma _{3}=\gamma }\left\Vert
\partial ^{\gamma _{1}}\eta _{k_{j}}(\xi )\cdot \partial ^{\gamma _{2}}|\xi
|^{-2}\cdot \partial ^{\gamma _{3}}(V_{\frac{n-2}{2}}(|\xi
|))^{-N}\right\Vert _{L^{2}}\left\Vert \Box _{k_{j}}\Delta
(A_{1})^{N}g_{k_{j}}\right\Vert _{L^{p}} \\
&\preceq &\sum\limits_{|\gamma |\leq \left[ n(\frac{1}{p}-\frac{1}{2})\right]
+1}\sum\limits_{\gamma _{1}+\gamma _{2}=\gamma }\left\Vert |\xi |^{-2}\cdot
\partial ^{\gamma _{1}}\eta _{k_{j}}(\xi )\cdot \partial ^{\gamma _{2}}(V_{%
\frac{n-2}{2}}(|\xi |))^{-N}\right\Vert _{L^{2}}\left\Vert \Box
_{k_{j}}\Delta (A_{1})^{N}g_{k_{j}}\right\Vert _{L^{p}} \\
&\preceq &\langle k_{j}\rangle ^{-(2-\frac{n-1}{2}N)}\Vert \Box
_{k_{j}}\Delta (A_{1})^{N}g_{k_{j}}\Vert _{L^{p}}
\end{eqnarray*}%
Combining above estimate with (\ref{eqo3.1}), Lemma \ref{lem3.2} is proved.

Next, we first verify the condition
\begin{equation*}
p_{1}\leq p_{2}
\end{equation*}%
and
\begin{equation*}
s_{1}\geq s_{2}+\sigma
\end{equation*}%
for $\forall q_{1},q_{2}\in \lbrack 1,+\infty )$. Let $f(x)$ be a nonzero
Schwartz function with $\mathrm{supp}\widehat{f}(\xi )\subset \{\xi :|\xi |<%
\frac{1}{2}\}$. Define
\begin{equation}\label{eqo3.5}
\widehat{f_{k_{j},\lambda }}(\xi )=\widehat{f}\left( \frac{\xi -k_{j}}{%
\lambda }\right)
\end{equation}%
for $\lambda \in (0,\rho ]$, where $\rho $ and $\{k_{j}\}$ are defined in
Lemma \ref{lem3.2}. By the definition of $f_{k_{j},\lambda }(x)$, we have
\begin{equation}\label{eqo3.20}
\Box _{k_{j}}f_{k_{j},\lambda }=f_{k_{j},\lambda }
\end{equation}%
and
\begin{equation}\label{eqo3.21}
\Box _{i}f_{k_{j},\lambda }(x)=0\quad \text{if}\quad i\neq k_{j}.
\end{equation}%
Then, by Lemma \ref{lem3.2}, we have
\begin{eqnarray*}
\left\Vert \Delta (A_{1})^{N}f_{k_{j},\lambda }\right\Vert
_{M_{p_{2},q_{2}}^{s_{2}}} &=&\langle k_{j}\rangle ^{s_{2}}\left\Vert \Box
_{k_{j}}\Delta (A_{1})^{N}f_{k_{j},\lambda }\right\Vert _{L^{p_{2}}} \\
&\simeq &\langle k_{j}\rangle ^{s_{2}+\sigma }\left\Vert f_{k_{j},\lambda
}\right\Vert _{L^{p_{2}}} \\
&\simeq &\langle k_{j}\rangle ^{s_{2}+\sigma }\lambda ^{n(1-\frac{1}{p_{2}}%
)}.
\end{eqnarray*}%
On the other hand
\begin{eqnarray*}
\Vert f_{k_{j},\lambda }\Vert _{M_{p_{1},q_{1}}^{s_{1}}} &=&\langle
k_{j}\rangle ^{s_{1}}\left\Vert \Box _{k_{j}}f_{k_{j},\lambda }\right\Vert
_{L^{p_{1}}} \\
&\simeq &\langle k_{j}\rangle ^{s_{1}}\Vert f_{k_{j},\lambda }\Vert
_{L^{p_{1}}} \\
&\simeq &\langle k_{j}\rangle ^{s_{1}}\lambda ^{n(1-\frac{1}{p_{1}})}.
\end{eqnarray*}%
By the assumption that $\Delta (A_{1})^{N}$ is bounded from $%
M_{p_{1},q_{1}}^{s_{1}}$ to $M_{p_{2},q_{2}}^{s_{2}}$, we have that
\begin{equation}
\langle k_{j}\rangle ^{s_{2}+\sigma }\lambda ^{n(1-\frac{1}{p_{2}})}\preceq
\langle k_{j}\rangle ^{s_{1}}\lambda ^{n(1-\frac{1}{p_{1}})}
\end{equation}%
for all $|k_{j}|$ sufficiently large and $0<\lambda \leq \rho $. Fix $k_{j}$
and let $\lambda $ go to 0. We have
\begin{equation*}
\lambda ^{n(1-\frac{1}{p_{2}})}\preceq \lambda ^{n(1-\frac{1}{p_{1}})},\quad
\text{for}\quad 0<\lambda \leq \rho .
\end{equation*}%
Thus, the condition $p_{1}\leq p_{2}$ must be hold. Moreover, when $\lambda $
is fixed and $k_{j}$ goes to infinity, we have
\begin{equation*}
\langle k_{j}\rangle ^{s_{2}+\sigma }\preceq \langle k_{j}\rangle
^{s_{1}},\quad \text{as}\quad k_{j}\rightarrow +\infty ,
\end{equation*}%
which yields $s_{2}+\sigma \leq s_{1}$.

For the condition of $q$, we first establish the following lemma.

\begin{lem}
\label{lem3.3}  For $j\in\mathbb{N}^{+}$, define
\begin{equation*}
\Lambda_{1,j}:=\{k\in\mathbb{Z}^{n}:|k|\in[j\pi+0.07,(j+1)\pi-0.07]\}
\end{equation*}
and
\begin{equation*}
\Lambda_{0,j}:=\{k\in\mathbb{Z}^{n}:|k|\in[j\pi,(j+1)\pi]\}.
\end{equation*}
When $j$ is big enough, we have
\begin{equation*}
|\Lambda_{1,j}|\geq C(n)|\Lambda_{0,j}|
\end{equation*}
where $C(n)$ is a positive constant depends only n.
\end{lem}

\textbf{Proof:} The proofs for $n=2$ and $n>2$ share the same idea. We prove
only the case $n=2$ explicitly and leave the proof of another case to the
reader.

By symmetry, we only need to consider the case $\{(x,y)\in \mathbb{R}%
^{2}:x,y\geq 0\}.$ For $j\in \mathbb{Z}^{+}$, we define
\begin{equation*}
I_{x>y}=\{(x,y)\in \mathbb{R}^{2}:|(x,y)|\in \lbrack j\pi ,j\pi
+0.07],x>y\geq 0\}
\end{equation*}%
\begin{equation*}
I_{x\leq y}=\{(x,y)\in \mathbb{R}^{2}:|(x,y)|\in \lbrack j\pi ,j\pi
+0.07],y\geq x\geq 0\}
\end{equation*}%
\begin{equation*}
II=\{(x,y)\in \mathbb{R}^{2}:|(x,y)|\in \lbrack j\pi +0.07,j\pi +\pi
-0.07],x,y\geq 0\}
\end{equation*}%
\begin{equation*}
III_{x>y}=\{(x,y)\in \mathbb{R}^{2}:|(x,y)|\in \lbrack j\pi +\pi -0.07,j\pi
+\pi ],x>y\geq 0\}
\end{equation*}%
\begin{equation*}
III_{x\leq y}=\{(x,y)\in \mathbb{R}^{2}:|(x,y)|\in \lbrack j\pi +\pi
-0.07,j\pi +\pi ],y\geq x\geq 0\}.
\end{equation*}%
Moreover, for $r,a>0$, $0\leq y\leq r$, we define an auxiliary function
\begin{equation*}
f_{r,a}(y)=\sqrt{(r+a)^{2}-y^{2}}-\sqrt{r^{2}-y^{2}}.
\end{equation*}%
Taking derivative we know that
\begin{equation*}
f_{r,a}^{\prime }(y)=y\left( \frac{1}{\sqrt{r^{2}-y^{2}}}-\frac{1}{\sqrt{%
(r+a)^{2}-y^{2}}}\right) \geq 0,
\end{equation*}%
and $f_{r,a}(y)$ is a monotone increasing function.

Then, for any $(x_{0},y_{0})\in I_{x>y}$, we have
\begin{equation*}
\left\vert \{x:(x,y_{0})\in I_{x>y}\}\right\vert =f_{j\pi ,0.07}(y_{0}).
\end{equation*}%
Therefore,
\begin{eqnarray*}
\max\limits_{(x_{0},y_{0})\in I_{x>y}}f_{j\pi ,0.07}(y_{0}) &=&f_{j\pi
,0.07}\left( \frac{1}{\sqrt{2}}(j\pi +0.07)\right)  \\
&=&\sqrt{(j\pi +0.07)^{2}-\frac{1}{\sqrt{2}}(j\pi +0.07)^{2}}-\sqrt{j\pi
^{2}-\frac{1}{\sqrt{2}}(j\pi +0.07)^{2}} \\
&=&\frac{0.14j+0.07^{2}}{\sqrt{(j\pi +0.07)^{2}-\frac{1}{\sqrt{2}}(j\pi
+0.07)^{2}}+\sqrt{j\pi ^{2}-\frac{1}{\sqrt{2}}(j\pi +0.07)^{2}}}.
\end{eqnarray*}%
It is obvious to see that $\lim\limits_{j\rightarrow +\infty
}\max\limits_{(x_{0},y_{0})\in I_{x>y}}f_{j\pi ,0.07}(y_{0})=\frac{0.14}{%
\sqrt{2}}<1$. Thus, for any $(x_{0},y_{0})\in I_{x>y}$ we have
\begin{equation*}
\left\vert \{x:(x,y_{0})\in I_{x>y}\}\right\vert \leq \max\limits_{(x,y)\in
I_{x>y}}f_{j\pi ,0.07}(y)<1,
\end{equation*}%
when $j$ is big enough.

On the other hand, for any $(x_{0},y_{0})\in II$, We have
\begin{equation*}
|\{x:(x,y_{0})\in II\}|=f_{j\pi +0.07,\pi -0.14}(y_{0}).
\end{equation*}%
By monotonicity of $f_{r,a}(y)$,
\begin{equation}
\min\limits_{(x_{0},y_{0})\in II}|\{x:(x,y_{0})\in II\}|=f_{j\pi +0.07,\pi -0.14}(0)=\pi -0.14>3.
\label{eqo3.4}
\end{equation}

Thus, for every $(x_{0},y_{0})\in I_{x>y}\bigcap \mathbb{Z}^{2}$, we have
\begin{equation*}
|\{(x,y_{0})\in \mathbb{Z}^{2}:(x,y_{0})\in I_{x>y}\}|=|\{(x_{0},y_{0})\}|=1
\end{equation*}%
and
\begin{equation*}
|\{(x,y_{0})\in \mathbb{Z}^{2}:(x,y_{0})\in II\}|\geq 3.
\end{equation*}%
Combing all above analysis, we have
\begin{equation*}
|\{(x,y)\in \mathbb{Z}^{2}:(x,y)\in II\}|\geq 3|\{(x,y)\in \mathbb{Z}%
^{2}:(x,y)\in I_{x>y}|.
\end{equation*}%
Now, we consider the domain $III_{x>y}$. By the same argument, for any $%
(x_{0},y_{0})\in III_{x>y}$, we have
\begin{eqnarray*}
|\{x:(x,y_{0})\in III_{x>y}\}| &\leq &\max\limits_{(x_{0},y_{0})\in
III_{x>y}}f_{j\pi +\pi -0.07,0.07}(y_{0}) \\
&=&f_{j\pi +\pi -0.07,0.07}\left( \frac{1}{\sqrt{2}}(j\pi +\pi )\right)  \\
&=&\sqrt{(j\pi +\pi )^{2}-\frac{1}{\sqrt{2}}(j\pi +\pi )^{2}}-\sqrt{(j\pi
+\pi -0.07)^{2}-\frac{1}{\sqrt{2}}(j\pi +\pi )^{2}} \\
&=&\frac{0.14j+0.07(2\pi -0.07)}{\sqrt{(j\pi +\pi )^{2}-\frac{1}{\sqrt{2}}%
(j\pi +\pi )^{2}}+\sqrt{(j\pi +\pi -0.07)^{2}-\frac{1}{\sqrt{2}}(j\pi +\pi
)^{2}}}.
\end{eqnarray*}%
It is easy to see
\begin{equation*}
\lim\limits_{j\rightarrow +\infty }\max\limits_{(x_{0},y_{0})\in
III_{x>y}}f_{j\pi +\pi -0.07,0.07}(y_{0})=\frac{0.14}{\sqrt{2}}<1.
\end{equation*}%
Thus, for any $(x_{0},y_{0})\in III_{x>y}$ we have
\begin{equation*}
\left\vert \{x:(x,y_{0})\in III_{x>y}\}\right\vert \leq
\max\limits_{(x,y)\in I_{x>y}}f_{j\pi +\pi -0.07,0.07}(y)<1,
\end{equation*}%
when $j$ is big enough. Moreover, it is obvious
\begin{equation*}
\frac{1}{\sqrt{2}}(j\pi +\pi )<j\pi +0.07
\end{equation*}%
when $j\geq 3$. So, for every $(x_{0},y_{0})\in III_{x>y}$,
\begin{equation*}
|\{x:(x,y_{0})\in II\}|=f_{j\pi +0.07,\pi -0.14}(y_{0})
\end{equation*}%
when $j\geq 3$. By (\ref{eqo3.4}), we can also obtain
\begin{equation*}
|\{(x,y)\in \mathbb{Z}^{2}:(x,y)\in II\}|\geq 3|\{(x,y)\in \mathbb{Z}%
^{2}:(x,y)\in I_{x>y}|.
\end{equation*}

On the other hand, for $I_{y\geq x}$ and $III_{y\geq x}$, by the same method
on the auxiliary function
\begin{equation*}
g_{r,a}(x)=\sqrt{(r+a)^{2}-x^{2}}-\sqrt{r^{2}-x^{2}},
\end{equation*}%
we can obtain that
\begin{equation*}
|\{(x,y)\in \mathbb{Z}^{2}:(x,y)\in II\}|\geq 3|\{(x,y)\in \mathbb{Z}%
^{2}:(x,y)\in I_{x\leq y}|
\end{equation*}%
and
\begin{equation*}
|\{(x,y)\in \mathbb{Z}^{2}:(x,y)\in II\}|\geq 3|\{(x,y)\in \mathbb{Z}%
^{2}:(x,y)\in III_{x\leq y}|.
\end{equation*}

Combing all above estimates, we have
\begin{equation*}
|\Lambda_{1,j}|\geq \frac{3}{7}|\Lambda_{0,j}|.
\end{equation*}

Next, we prove the necessary conditions when $q_{1}\leq q_{2}$. By (2.1) and
above analysis, the boundedness of $\Delta (A_{1})^{N}$ from $%
M_{p_{1},q_{1}}^{s_{1}}$ to $M_{p_{2},q_{2}}^{s_{2}}$ must hold for $%
p_{1}\leq p_{2},s_{2}+\sigma \leq s_{1}$ and $q_{1}\leq q_{2}$. Therefore,
we only need to consider the case $q_{1}>q_{2}$. Let $M$ be a large positive
number. Define
\begin{equation*}
F_{M}(x)=\sum\limits_{100<|k_{j}|<M}a_{j}f_{k_{j},\rho }(x)
\end{equation*}%
where $a_{j}>0$ are constants to be chosen later and $f_{k_{j},\rho }(x)$
are defined in (\ref{eqo3.5}) with all $k_{j}$ satisfy
\begin{equation*}
|k_{j}|\in \lbrack L\pi +0.07,L\pi +\pi -0.07]
\end{equation*}%
for some $L\in \mathbb{N^{+}}.$

By (\ref{eqo3.20}) (\ref{eqo3.21}) and the almost orthogonality of $\{\varphi_{k}\}$, we have
\begin{eqnarray*}
\left\Vert \Delta (A_{1})^{N}F_{M}\right\Vert _{M_{p_{2},q_{2}}^{s_{2}}}
&=&\left( \sum\limits_{k\in \mathbb{Z}^{n}}a_{j}^{q_{2}}\langle k\rangle
^{s_{2}q_{2}}\left\Vert \Box _{k}\Delta (A_{1})^{N}F_{M}\right\Vert
_{L^{p_{2}}}^{q_{2}}\right) ^{\frac{1}{q_{2}}} \\
&\simeq &\left( \sum\limits_{100<|k_{j}|<M}a_{j}^{q_{2}}\langle k_{j}\rangle
^{s_{2}q_{2}}\langle k_{j}\rangle ^{\sigma q_{2}}\left\Vert f_{k_{j},\rho
}(x)\right\Vert _{L^{p_{2}}}^{q_{2}}\right) ^{\frac{1}{q_{2}}} \\
&\simeq &\left( \sum\limits_{100<|k_{j}|<M}a_{j}^{q_{2}}\langle k_{j}\rangle
^{s_{2}q_{2}+\sigma q_{2}}\right) ^{\frac{1}{q_{2}}}
\end{eqnarray*}%
and
\begin{eqnarray*}
\left\Vert F_{M}\right\Vert _{M_{p_{1},q_{1}}^{s_{1}}} &=&\left(
\sum\limits_{k\in \mathbb{Z}^{n}}a_{j}^{q_{1}}\langle k\rangle
^{s_{1}q_{1}}\Vert \Box _{k}F_{M}\Vert _{L^{p_{1}}}^{q_{1}}\right) ^{\frac{1%
}{q_{1}}} \\
&\simeq &\left( \sum\limits_{100<|k_{j}|<M}a_{j}^{q_{1}}\langle k_{j}\rangle
^{s_{1}q_{1}}\left\Vert f_{j,\rho }(x)\right\Vert
_{L^{p_{1}}}^{q_{1}}\right) ^{\frac{1}{q_{1}}} \\
&\simeq &\left( \sum\limits_{100<|k_{j}|<M}a_{j}^{q_{1}}\langle k_{j}\rangle
^{s_{1}q_{1}}\right) ^{\frac{1}{q_{1}}}.
\end{eqnarray*}%
By the assumption that $\Delta (A_{1})^{N}$ is bounded form $%
M_{p_{1},q_{1}}^{s_{1}}$ to $M_{p_{2},q_{2}}^{s_{2}}$, we have
\begin{equation}
\left( \sum\limits_{100<|k_{j}|<M}a_{j}^{q_{2}}\langle k_{j}\rangle
^{s_{2}q_{2}+\sigma q_{2}}\right) ^{\frac{1}{q_{2}}}\preceq \left(
\sum\limits_{100<|k_{j}|<M}a_{j}^{q_{1}}\langle k_{j}\rangle
^{s_{1}q_{1}}\right) ^{\frac{1}{q_{1}}}.
\end{equation}%
By choosing $a_{j}=\langle k_{j}\rangle ^{\frac{s_{1}q_{1}-(s_{2}+\sigma
)q_{2}}{q_{1}-q_{2}}}$, we obtain
\begin{equation}
\left( \sum\limits_{100<|k_{j}|<M}\langle k_{j}\rangle ^{\frac{\lbrack
s_{1}-(s_{2}+\sigma )]q_{1}q_{2}}{q_{2}-q_{1}}}\right) ^{\frac{1}{q_{2}}%
}\preceq \left( \sum\limits_{100<|k_{j}|<M}\langle k_{j}\rangle ^{\frac{%
\lbrack s_{1}-(s_{2}+\sigma )]q_{1}q_{2}}{q_{2}-q_{1}}}\right) ^{\frac{1}{%
q_{1}}}.
\end{equation}%
By the assumption $q_{1}>q_{2}$, the above series converges as $M\rightarrow
+\infty $. By Lemma \ref{lem3.3}, we have
\begin{equation}
\sum\limits_{100<|k_{j}|<M}\langle k_{j}\rangle ^{\frac{\lbrack
s_{1}-(s_{2}+\sigma )]q_{1}q_{2}}{q_{2}-q_{1}}}\simeq
\sum\limits_{100<|k|<M}\langle k\rangle ^{\frac{\lbrack s_{1}-(s_{2}+\sigma
)]q_{1}q_{2}}{q_{2}-q_{1}}}.
\end{equation}%
Therefore, it must yield
\begin{equation}
\frac{\lbrack s_{1}-(s_{2}+\sigma )]q_{1}q_{2}}{q_{2}-q_{1}}<-n,
\end{equation}%
which is equivalent to $s_{1}+\frac{n}{q_{1}}>s_{2}+\sigma +\frac{n}{q_{2}}$%
. Theorem 1.1 is proved.

\end{document}